\documentclass[12pt,reqno]{amsart}
\usepackage[T1]{fontenc}
\usepackage[utf8]{inputenc}
\usepackage{lmodern,amsmath,amssymb,mathtools,microtype}
\usepackage{xurl}
\usepackage[hidelinks]{hyperref}
\newtheorem{theorem}{Theorem}[section]
\newtheorem{corollary}[theorem]{Corollary}

\theoremstyle{remark}

\newcommand{\C}{\mathbb{C}}

\newcommand{\eps}{\varepsilon}
\newcommand{\norm}[1]{\lVert #1\rVert}

\title[A strict numerical-index bound]{A strict numerical-index bound for complex Lebesgue spaces}
\author{M.A.~Perelmuter}
\address{SCAD Soft Ltd., 3a Osvity Street, Kyiv, 03037, Ukraine}
\email{mikeperelmuter@gmail.com}
\date{}
\subjclass[2020]{Primary 47A12; Secondary 46B20, 46E30}
\keywords{Numerical index, numerical radius, $L_p$-spaces}

\begin{document}
\begin{abstract}
We prove that the numerical index of a complex Lebesgue space $L_p(\mu)$ of dimension at least two is strictly less than $p^{-1/p}q^{-1/q}$ for every $1<p<\infty$, $p\ne 2$, where $q$ is the conjugate exponent of $p$. This disproves, in the complex case, a conjecture of Mart\'in, Mer\'i and Popov. The proof is based on a perturbation of a rank-one operator on the two-dimensional complex space $\ell_p^2$.
\end{abstract}
\maketitle

\section{Introduction}

The numerical index of a Banach space is a classical isometric invariant relating the numerical radius and the operator norm. Given a real or complex Banach space $X$, we denote by $X^*$ its dual space and by $\mathcal L(X)$ the space of all bounded linear operators on $X$. The numerical radius of $T\in\mathcal L(X)$ is defined by
$$v(T)=\sup\bigl\{|x^*(Tx)|:\ x\in X,\ x^*\in X^*,\ \|x\|=\|x^*\|=x^*(x)=1\bigr\}.$$
The numerical index of $X$ is given by
$$n(X)=\inf\bigl\{v(T):\ T\in\mathcal L(X),\ \|T\|=1\bigr\}.$$
Equivalently, $n(X)$ is the greatest constant $c\ge0$ such that
$$c\|T\|\le v(T),\qquad T\in\mathcal L(X).$$
For background on the numerical index and related open problems, we refer to the survey \cite{KMP}.

The computation of the numerical index of Lebesgue spaces $L_p(\mu)$ for $1<p<\infty$, $p\ne2$, is one of the longstanding problems in the theory of numerical ranges. Despite considerable progress (see \cite{EdDariKhamsi,MM,MMP,Chiclana} and the references therein), the exact value of the numerical index of infinite-dimensional Lebesgue spaces remains unknown in both the real and the complex case.

Let $q=p/(p-1)$ be the conjugate exponent of $p$ and set
$$M_p=\max_{0\le t\le1}\frac{|t^{p-1}-t|}{1+t^p},\qquad\kappa_p=p^{-1/p}q^{-1/q}.$$
Mart\'in, Mer\'i and Popov \cite[Problem~5.1]{MMP} conjectured that, for infinite-dimensional Lebesgue spaces,
$$n(L_p(\mu))=\begin{cases}
M_p & \text{in the real case},\\
\kappa_p & \text{in the complex case}.
\end{cases}$$

In a recent preprint, Chiclana \cite{Chiclana} proves that
$$n(\ell_p^2(\mathbb R))=M_p$$
for every $1<p<\infty$. However, in a recent manuscript, Colbrook \cite{Colbrook} shows that
$$n(\ell_p^3(\mathbb R))<M_p$$
for every $1<p<\infty$, $p\ne2$, thereby disproving the conjecture in the real case, since $n(L_p(\mu))\le n(\ell_p^3)$ for infinite-dimensional real $L_p(\mu)$ (see \cite[Proposition~6.1(a),(b)]{MMPR}).

The purpose of this note is to disprove the conjecture in the complex case. More precisely, we prove that
$$n(\ell_p^2(\mathbb C))<\kappa_p,\qquad 1<p<\infty,\quad p\ne2.$$
It follows that
$$n(L_p(\mu))<\kappa_p$$
for every complex Lebesgue space $L_p(\mu)$ of dimension at least two.

Our result also answers negatively the question whether $n(\ell_p^2(\C))=\kappa_p$, posed as Problem~3 in the survey \cite{KMP}.

Our proof uses the perturbation mechanism employed by Colbrook in the real case. For $p>2$, we construct a perturbation of a norm-one operator whose numerical radius increases at most quadratically in the perturbation parameter, whereas its norm increases at least at a rate of order $\varepsilon^q$, with $q<2$. In contrast to the real case, a two-dimensional space suffices.

Throughout the paper, $(\Omega,\Sigma,\mu)$ denotes a measure space, and $L_p(\mu)$ and $\ell_p^m$ denote complex Banach spaces unless otherwise stated. For $x\in L_p(\mu)$ we write
$$x^\#=|x|^{p-1}\operatorname{sign}(\overline{x}).$$
The numerical radius of an operator $T\in\mathcal L(L_p(\mu))$ can then be expressed as
$$v(T)=\sup_{\|x\|_p=1}\left|\int\limits_\Omega x^\# Tx\,\mathrm d\mu\right|.$$

\section{Main result}

\begin{theorem}\label{thm:main}
Let $1<p<\infty$, $p\neq2$. Then $n(\ell_p^2)<\kappa_p$.
\end{theorem}
\begin{proof}
The numerical radii of an operator and its adjoint coincide, as do their operator norms. In finite dimension, every operator on the dual space is an adjoint. Thus a finite-dimensional Banach space and its dual have the same numerical index. Here $(\ell_p^2)^*=\ell_q^2$. Since $n(\ell_p^2)=n(\ell_q^2)$ and $\kappa_p=\kappa_q$, it suffices to consider the case $p>2$, so that $1<q<2$.

We perturb the operator $\left(\begin{smallmatrix}0&0\\1&0\end{smallmatrix}\right)$, whose norm is $1$ and whose numerical radius is $\kappa_p$ (see \cite[Lemma~2]{MM}). For $\eps>0$ define
$$T_\eps=\begin{pmatrix}\eps(p-1)&0\\1&-\eps\end{pmatrix}.$$

\emph{Estimate of the numerical radius.}
Take $x=(x_1,x_2)$ with $\norm{x}_p=1$ and put $\lambda=|x_1|^p\in[0,1]$, so that $|x_2|^p=1-\lambda$. Since $x_j^\#x_j=|x_j|^p$, we have
$$x^\#T_\eps x=\eps(p-1)x_1^\#x_1+x_2^\#x_1-\eps x_2^\#x_2=x_2^\#x_1+\eps\bigl((p-1)\lambda-(1-\lambda)\bigr)
=x_2^\#x_1+\eps(p\lambda-1).$$
Moreover, since $(p-1)/p=1/q$,
$$|x_2^\#x_1|=|x_2|^{p-1}|x_1|=\lambda^{1/p}(1-\lambda)^{1/q}=:f(\lambda).$$
The triangle inequality gives
$$v(T_\eps)\le\max_{0\le\lambda\le1}\bigl(f(\lambda)+\eps|p\lambda-1|\bigr).$$

Let $a=1/p$ and $b=1/q$, so $a+b=1$. Differentiating on $(0,1)$ gives
$$f'(\lambda)=a\lambda^{a-1}(1-\lambda)^b-b\lambda^a(1-\lambda)^{b-1} = \lambda^{a-1}(1-\lambda)^{b-1}(a-\lambda).$$
Its maximum is $f(a)=a^ab^b=\kappa_p$. Differentiating once more,
\begin{align*}
f''(\lambda) =&a(a-1)\lambda^{a-2}(1-\lambda)^b -2ab\lambda^{a-1}(1-\lambda)^{b-1}+b(b-1)\lambda^a(1-\lambda)^{b-2}\\
=&-ab\lambda^{a-2}(1-\lambda)^{b-2}\bigl((1-\lambda)^2+2\lambda(1-\lambda)+\lambda^2\bigr)\\
=&-\frac1{pq}\lambda^{1/p-2}(1-\lambda)^{1/q-2} \le-\frac1{pq}.
\end{align*}
The inequality follows because both exponents are negative and $\lambda,1-\lambda\in(0,1)$. By Taylor's formula, for $\lambda\in(0,1)$ there is a point $\xi$ between $a$ and $\lambda$ such that, with $h=\lambda-a$,
$$f(\lambda)=f(a)+f'(a)h+\tfrac12 f''(\xi)h^2\le\kappa_p-\frac{h^2}{2pq}.$$
The same bound holds at the endpoints by continuity. Since $p\lambda-1=ph$, we obtain
$$f(\lambda)+\eps|p\lambda-1|\le\kappa_p-\frac{h^2}{2pq}+p\eps|h|=\kappa_p+\frac{p^3q}{2}\eps^2
-\frac{(|h|-p^2q\eps)^2}{2pq}\le\kappa_p+\frac{p^3q}{2}\eps^2.$$
Consequently,
$$v(T_\eps)\le\kappa_p+\frac{p^3q}{2}\eps^2.$$

\emph{Estimate of the operator norm.}
Let $t=\eps^{q-1}$ and $x=(1,-t)$. The identities $p(q-1)=q$ and $1+(q-1)=q$ yield $t^p=\eps^q$, $\eps t=\eps^q$. Hence $T_\eps x=(\eps(p-1),1+\eps^q)$ and
$$\norm{x}_p = (1+\eps^q)^{1/p},\qquad \norm{T_\eps x}_p \ge 1+\eps^q. $$
By the definition of the operator norm,
$$\norm{T_\eps}\ge\frac{\norm{T_\eps x}_p}{\norm{x}_p} \ge(1+\eps^q)^{1-1/p}=(1+\eps^q)^{1/q}.$$

\emph{Comparison of the estimates.}
Put $A=p^3q/2$. Combining the
above estimates, we obtain
$$
n(\ell_p^2)\le\frac{v(T_\eps)}{\norm{T_\eps}}
\le\frac{\kappa_p+A\eps^2}{(1+\eps^q)^{1/q}}.
$$
By concavity,
$$
(1+s)^{1/q}\ge 1+(2^{1/q}-1)s,\qquad s\in[0,1].
$$
Hence, for $\eps\in(0,1)$,
$$
n(\ell_p^2)
\le\frac{\kappa_p+A\eps^2}{(1+\eps^q)^{1/q}}
\le\frac{\kappa_p+A\eps^2}
{1+(2^{1/q}-1)\eps^q}.
$$
Since $1<q<2$, we can take
$$
\eps=\left(\frac{\kappa_p(2^{1/q}-1)}{2A}\right)^{1/(2-q)}
\in(0,1).
$$
For this choice of $\eps$,
$$
n(\ell_p^2)
\le \kappa_p\frac{1+\frac12(2^{1/q}-1)\eps^q}
{1+(2^{1/q}-1)\eps^q}
<\kappa_p.
$$
\end{proof}

\begin{corollary}
$n(L_p(\mu))<\kappa_p$ whenever $\dim L_p(\mu)\ge2$, $1<p<\infty$, and $p\ne2$.
\end{corollary}

\begin{proof}
If $L_p(\mu)$ has finite dimension $m\ge2$, it is isometric to $\ell_p^m$, and $n(\ell_p^m)\le n(\ell_p^2)$ by \cite[Proposition~6.1(a)]{MMPR}.
If it is infinite-dimensional, then $n(L_p(\mu))=\inf_m n(\ell_p^m)$ by \cite[Proposition~6.1(a),(b)]{MMPR}.
\end{proof}

\section*{Acknowledgments}
The author used ChatGPT (OpenAI) and Claude (Anthropic) to assist in
finding suitable parameters for the matrix $T_\eps$, refining the
mathematical arguments, and improving their presentation.
The author has verified all mathematical results and proofs and takes
full responsibility for the content of the paper.

\end{document}